\documentclass[12pt]{article}
\usepackage{mathrsfs}
\usepackage{amsfonts,amsmath,amsthm, amssymb}
\usepackage{latexsym, euscript, epic, eepic}
\usepackage{lineno}

\newtheorem{cor}{Corollary}[section]
\newtheorem{defn}{Definition}[section]
\newtheorem{prop}{Proposition}[section]
\newtheorem{thm}{Theorem}[section]

\begin{document}

\title{Induced  topological pressure for topological dynamical
systems
  }
\author{Zhitao Xing$^{1,2},$ Ercai Chen$^{1,3}.$ \\
   \small   1 School of Mathematical Science, Nanjing Normal University,\\
    \small   Nanjing 210023, Jiangsu, P.R. China\\
     \small    e-mail: ecchen@njnu.edu.cn \\
     \small    e-mail: xzt-303@163.com \\
      \small  2 School of Mathematics and  Statistics,
                 Zhaoqing University,\\
       \small   Zhaoqing 526061, Guangdong, P.R. China\\
        \small 3 Center of Nonlinear Science, Nanjing University,\\
         \small   Nanjing 210093, Jiangsu, P.R. China\\}
\date{}
\maketitle

\begin{center}
 \begin{minipage}{120mm}
{\small {\bf Abstract.} In this paper, inspired  by  the  article
[5], we introduce the induced topological pressure for a topological
dynamical system. In particular, we prove a variational principle
for the induced topological pressure. }
\end{minipage}
 \end{center}

\vskip0.5cm
\section{INTRODUCTION AND  MAIN RESULT  }

\quad \quad The present paper is devoted to the study of  the
induced topological pressure for topological dynamical systems.
Before stating our main result, we first give some notation  and
background about the  induced topological pressure. By a
topological dynamical system (TDS) $(X,f)$, we mean a compact
metric space $(X,d)$ together with a continuous map
$f:X\rightarrow X.$ Recall that $C(X,\mathbb{R})$
  is the Banach algebra of real-valued continuous functions of $X$
equipped with the supremum norm. For $\varphi \in C(X,\mathbb{R}),
n\geq 1$, let $(S_{n}\varphi)(x):=\sum
\limits_{i=0}^{n-1}\varphi(f^{i}x)$ and for $\psi \in
C(X,\mathbb{R})$ with $\psi >0$, let $m:=\min\{\psi(x): x\in X\}$.
We denote by $ M(X,f)$ all $f$-invariant Borel probability
measures on $X$ endowed with the weak-star topology.

  Topological pressure is a  basic
notion of the thermodynamic formalism. It first introduced by
Ruelle [11] for expansive topological  dynamical systems, and  later
by Walters [1,9,10] for the general case.
 The variational principle  established by Walters can be stated  as follows: Let $(X,f)$ be a TDS, and let $\varphi \in
C(X,\mathbb{R})$, $ P(\varphi)$ denote the topological pressure of
$\varphi.$ Then
\begin{equation}\label{tag-1}
P(\varphi)=\sup\{h_{\mu}(f)+\int \varphi d\mu: \mu \in M(X,f)\}.
\end{equation}
where $h_{\mu}(f)$ denotes the  measure-theoretical entropy of
$\mu.$    The theory of topological pressure and its variational
principle  plays a fundamental role in statistics, ergodic theory,
and the theory of dynamical systems [3,9,13]. Since the works of Bowen [4] and
Ruelle [12], the topological pressure has become a basic tool in
the dimension theory of dynamical systems [8,14].

 Recently  Jaerish, Kesseb\"{o}hmer and
Lamei [5] introduced  the  notion of the induced topological pressure of a countable
Markov shift, and established a variational principle for it.
One important feature of this pressure is the freedom in choosing
a scaling function, and this is applied to large deviation theory
and fractal geometry. In this paper we present the induced topological
pressure for a topological dynamical system and consider the
relation between   it and the topological
pressure. We set up a variational principle for the induced
topological pressure. As an application, we
will point out that the BS dimension is a  special case of the induced
topological pressure.

  Let  $(X,f)$ be a TDS. For $n\in \mathbb{N}$, the $n$th Bowen
metric $d_{n}$ on $X$ is defined by
$$d_{n}(x,y)=\max \{d(f^{i}(x),f^{i}(y)): i=0,1,\ldots, n-1 \}.$$
For every $\epsilon
>0$, we denote by $B_{n}(x,\epsilon),\overline{B}_{n}(x,\epsilon) $ the open (resp. closed) ball of
radius $\epsilon$ and order $n$ in the metric $d_{n}$ around $x$,
i.e., $$B_{n}(x,\epsilon)= \{y\in X : d_{n}(x,y)<\epsilon\} \text{
and } \overline{B}_{n}(x,\epsilon)= \{y\in X : d_{n}(x,y)\leq
\epsilon\}. $$ Let $Z\subseteq X$ be a non-empty set. A subset
$F_{n}\subset X$ is called an $(n, \epsilon)$-spanning set of $Z$
if for any $y\in Z$, there exists $x \in F_{n} $ with
$d_{n}(x,y)\leq \epsilon$.\ A subset $E_{n}\subset Z$ is called an
$(n,\epsilon)$-separated set of $Z$ if $x,y\in E_{n}, x\neq y$
implies $d_{n}(x,y)>\epsilon$.

 Now we define a new notion, the \textit{induced
topological  pressure} which extends the definition in [5] for
topological Markov shifts if the Markov shift is compact, as
follows.
\begin{defn}
Let $(X,f)$ be a TDS and $ \varphi,\psi \in C(X,\mathbb{R})$ with
$\psi>0$.  For $ T>0$,  define
$$S_{T}=\{n\in \mathbb{N}: \exists x\in X \text { such  that  }
S_{n}\psi(x)\leq T \text{ and }S_{n+1}\psi(x)>T\}.$$ For $n\in
S_{T}$, define
$$X_{n}=\{x\in X: S_{n}\psi(x)\leq T \text{ and }S_{n+1}\psi(x)>T \}.$$
Let
$$Q_{\psi ,T}(f,\varphi, \epsilon)=
\inf\left\{\sum\limits_{n\in S_{T}}\sum \limits_{x\in F_{n}}\exp
(S_{n}\varphi)(x): F_{n} \ is \ an \ (n,\epsilon)\text{-spanning
set of } X_{n},n\in S_{T} \right\}.$$  We define the
$\psi$-induced topological pressure of $\varphi $ (with respect to
$f$) by
\begin{equation}\label{tag-1}
P_{\psi}(\varphi)=\lim \limits_{\epsilon\rightarrow 0}\limsup
\limits_{T \rightarrow \infty}\frac{1}{T}\log Q_{\psi
,T}(f,\varphi, \epsilon)
\end{equation}
\end {defn}
\noindent {\bf Remarks.}\\  ($\romannumeral1$) Let
  $[\frac{T}{m}]$ denote the integer part of $\frac{T}{m}$. Then for $n\in S_{T}$, $n\leq [\frac{T}{m}]+1$, i.e., $S_T$ is a finite set.\\($\romannumeral2$) If
  $0<\epsilon_{1}<\epsilon_{2}$,
then $Q_{\psi ,T}(f,\varphi, \epsilon_{1})\geq Q_{\psi
,T}(f,\varphi, \epsilon_{2})$, which implies the existence of the
limit in (1.2) and  $P_{\psi}(\varphi)> -\infty$.\\
($\romannumeral3$) $P_1(\varphi)=P(\varphi)$.

  The
 variational principle for  induced topological pressure is stated as
 follows.
\begin{thm}
Let $(X,f)$ be a TDS and   $\varphi, \psi \in C(X,\mathbb{R})$ with
$\psi> 0$. Then
\begin{equation}\label{tag-1}
P_{\psi}(\varphi)=\sup\left\{\frac{h_{\nu}(f)}{\int \psi
d\nu}+\frac{\int \varphi d\nu}{\int \psi d\nu}: \nu \in M(X ,f
)\right\}.
\end{equation}
\end{thm}

This  paper is organized  as follows. In Section 2, we provide an
equivalent definition of induced topological pressure. We prove Theorem 1.1 in
Section 3. We point out  that the BS dimension is a special case of
the induced topological pressure in Section 4. In Section 5, we study the equilibrium
measures for  the induced topological pressure.
\section{AN EQUIVALENT DEFINITION  }

\quad\quad In this section,  we obtain  an equivalent definition
of  the induced topological pressure by using separated sets (from now on, we omit the word `topological' if no confusion can arise).
\begin{prop}
Let $(X,f)$ be a TDS and $\varphi, \psi \in C(X,\mathbb{R})$ with
$\psi>0$.  For $T>0$,  define $$P_{\psi ,T}(f,\varphi, \epsilon)=
\sup\left\{\sum\limits_{n\in S_{T}}\sum \limits_{x\in E_{n}}\exp
(S_{n}\varphi)(x): E_{n} \ is \ an \ (n,\epsilon)\text{-separated
set of } X_n,  n\in S_T \right\}.$$  Then
\begin{equation}\label{tag-1}P_{\psi}(\varphi)=\lim \limits_{\epsilon\rightarrow 0}\limsup
\limits_{T \rightarrow \infty}\frac{1}{T}\log P_{\psi
,T}(f,\varphi, \epsilon)
\end{equation}
\end{prop}
\noindent\textbf{Proof.} We note that since  the map
$\epsilon\mapsto \limsup \limits_{T \rightarrow
\infty}\frac{1}{T}\log P_{\psi ,T}(f,\varphi, \epsilon)$ is
nondecreasing, the limit in (2.4) is well defined when
$\epsilon\rightarrow 0$. For $n\in S_{T}$, let $E_{n}$ be an
$(n,\epsilon)$-separated set of $X_{n}$ which fails to be
$(n,\epsilon)$-separated  when any point of $X_{n}$ is added. Then
$E_{n}$ is an $(n,\epsilon)$-spanning set of $X_{n}$. Therefore
$$Q_{\psi ,T}(f,\varphi, \epsilon)\leq P_{\psi ,T}(f,\varphi,
\epsilon)$$ and $$P_{\psi}(\varphi)\leq \lim
\limits_{\epsilon\rightarrow 0}\limsup \limits_{T \rightarrow
\infty}\frac{1}{T}\log P_{\psi ,T}(f,\varphi, \epsilon).$$

   To show the reverse
inequality, for any $\epsilon>0$, we choose $\delta>0$ small
enough so that \begin{flalign} d(x,y)\leq \frac{\delta}{2}
\Rightarrow |\varphi(x)-\varphi(y)|<\epsilon.
\end{flalign}
For $n\in S_{T}$, let  $E_{n}$ be an $(n,\delta)$-separated set of
$X_{n}$ and $F_{n}$ an $(n, \frac{\delta}{2})$-spanning set of
$X_{n}$. Define $\phi:E_{n}\rightarrow F_{n}$ by choosing, for
each $x\in E_{n}$, some point $\phi(x)\in F_{n}$ with
$d_{n}(\phi(x),x)\leq \frac{\delta}{2}$. Then $\phi$ is injective.
\\ Therefore,
\begin{flalign*} &\sum\limits_{n\in S_{T}}\sum\limits_{y\in
F_{n}}\exp(S_{n}\varphi)(y)\\
\geq&\sum\limits_{n\in S_{T}}\sum\limits_{y\in\phi E_{n}}\exp(S_{n}\varphi)(y)\\
\geq&\sum\limits_{n\in S_{T}}(\min\limits_{x\in
E_{n}}\exp((S_{n}\varphi)(\phi
x)-(S_{n}\varphi)(x)))\sum\limits_{x\in
E_{n}}\exp(S_{n}\varphi)(x)
\\\geq& \exp(-(\frac{T}{m}+1)\epsilon)\sum\limits_{n\in S_{T}}\sum\limits_{x\in E_{n}}\exp(S_{n}\varphi)(x).
\end{flalign*}
  We conclude that
$$ \lim \limits_{\delta\rightarrow
0}\limsup \limits_{T \rightarrow \infty}\frac{1}{T}\log Q_{\psi
,T}(f,\varphi, \frac{\delta}{2})\geq -\frac{1}{m}\epsilon+\lim
\limits_{\delta\rightarrow 0}\limsup \limits_{T \rightarrow
\infty}\frac{1}{T}\log P_{\psi ,T}(f,\varphi, \delta).$$ As
$\epsilon\rightarrow 0 $, we have
$$P_{\psi}(\varphi)\geq \lim \limits_{\delta\rightarrow
0}\limsup \limits_{T \rightarrow \infty}\frac{1}{T}\log P_{\psi
,T}(f,\varphi, \delta).$$
\section{ THE PROOF OF THEOREM 1.1}

\quad\quad In this section, we give the proof of Theorem 1.1.
Firstly, we study   the relation between $P_{\psi}(\varphi)$ and
$P(\varphi)$, which will be needed for the proof of Theorem 1.1.
The following Theorem 3.1 is very similar to Theorem 2.1 of [5],
and  it is a generalization of this theorem in the case of a
compact topological Markov shift.
\begin{thm}
Let $(X,f)$ be a TDS and  $\varphi, \psi \in C(X,\mathbb{R})$ with
$\psi> 0$.  For $T>0$,  define $$G_{T}=\{n\in \mathbb{N}: \exists
x\in X \text { such  that  } S_{n}\psi(x)> T\}.$$ For $n\in
G_{T}$, define $$Y_{n}=\{x\in X: S_{n}\psi(x)>T\}.$$ Let
$$R_{\psi ,T}(f,\varphi, \epsilon)=
\sup\left\{\sum\limits_{n\in G_{T}}\sum \limits_{x\in
E^{'}_{n}}\exp (S_{n}\varphi)(x): E^{'}_{n} \text { is \ an
}(n,\epsilon)\text{-separated  set of } Y_{n}, n\in G_T \right
\}.$$  We have
\begin{flalign}
P_{\psi}(\varphi)=\inf\{\beta \in \mathbb{R}: \lim
\limits_{\epsilon\rightarrow 0}\limsup \limits_{T \rightarrow
\infty}R_{\psi ,T}(f,\varphi-\beta\psi, \epsilon)<\infty\}.
\end{flalign}
Here we make the convention that $\inf \emptyset =\infty$.
\end{thm}
\noindent\textbf{Proof.}   For $n\in \mathbb{N},x\in X$, we define
$m_n(x)$ to be the unique positive integer
 such that
$$(m_n(x)-1)\|\psi\|<S_{n}\psi(x)\leq m_n(x)\|\psi\|.$$  Observing that $$ \exp(-\beta
\|\psi\| m_n(x))\exp(-|\beta|\|\psi\|)\leq \exp(-\beta
S_{n}\psi(x))\leq \exp(-\beta \|\psi\|
m_n(x))\exp(|\beta|\|\psi\|)$$ for all $ x\in X$. For
$\xi_T=\{\xi_n: X\to\mathbb{R}\}_{n\in G_T}$, we define
\begin{flalign*} &R_{\psi ,T}(f,\varphi, \xi_T,\epsilon)\\=&
\sup\left\{\sum\limits_{n\in G_{T}}\sum \limits_{x\in
E^{'}_{n}}\exp ((S_{n}\varphi)(x)-\xi_n(x)): E^{'}_{n} \text { is
\ an } (n,\epsilon)\text{-separated  set of } Y_{n}, n\in G_T
\right\}.\end{flalign*}We conclude that
$$\lim \limits_{\epsilon\rightarrow 0}\limsup \limits_{T
\rightarrow \infty}R_{\psi ,T}(f,\varphi-\beta \psi,
\epsilon)<\infty $$ if and only if
$$\lim \limits_{\epsilon\rightarrow 0}\limsup \limits_{T
\rightarrow \infty}R_{\psi ,T}(f,\varphi,\{-\beta \|\psi\|
m_n\}_{n\in G_T}, \epsilon)<\infty.$$
 Hence, it will be sufficient to verify that
$$P_{\psi}(\varphi)=\inf\{\beta \in \mathbb{R}: \lim \limits_{\epsilon\rightarrow 0}\limsup \limits_{T
\rightarrow \infty}R_{\psi ,T}(f,\varphi,\{-\beta \|\psi\|
m_n\}_{n\in G_{T}}, \epsilon)<\infty\}.$$ By the equivalent
definition of $P_{\psi}(\varphi),$ for every $ \delta>0 , \beta\in
\mathbb{R}$ with $\beta<P_{\psi}(\varphi)-\delta$, there exists an
$\epsilon_{0}>0$ with   $$\beta+\delta<\limsup \limits_{T
\rightarrow \infty}\frac{1}{T}\log P_{\psi ,T}(f,\varphi,
\epsilon)\leq P_{\psi}(\varphi), \ \ \forall \epsilon \in (0,
\epsilon_{0}),$$ and we can find  a sequence $\{T_{j}\}_{j\in
\mathbb{N}}$ such that  for every $j\in \mathbb{N},$
$T_{j+1}-T_{j}>2\|\psi\|$ and for each $j\in \mathbb{N}$, there
exists an $E_{T_{j}}=\bigcup \limits_{n\in S_{T_{j}}}E_{n}$ with
$$\sum\limits_{n\in S_{T_{j}}}\sum
\limits_{x\in E_{n}}\exp (S_{n}\varphi)(x)\geq
\exp(T_{j}(\beta+\frac{\delta}{2})).$$  Since for $j\in
\mathbb{N}, n\in S_{T_{j}}, x\in E_{n}$,
$T_j-\|\psi\|<S_{n}\psi(x)\leq T_j$, we have
$$S_{T_i}\cap S_{T_j}=\emptyset,i\neq j$$
and
$$|\|\psi\| m_n(x)-T_{j}|<2\|\psi\|.$$
It follows that
\begin{flalign*}
&R_{\psi ,T}(f,\varphi,\{-\beta
\|\psi\| m_n\}_{n\in G_T},\epsilon)\\
\geq& \sum\limits_{j\in \mathbb{ N },\ T_{j}-\|\psi\|>T}\sum
\limits_{n\in S_{T_{j}}}\sum
\limits_{x\in E_{n}}\exp ((S_{n}\varphi)(x)-\beta \|\psi\| m_n(x))\\
\geq& \exp(-2|\beta|\|\psi\|)\sum\limits_{j\in \mathbb{N},\
T_{j}-\|\psi\|> T}\sum\limits_{n\in
S_{T_{j}}}\sum \limits_{x\in E_{n}}\exp ((S_{n}\varphi)(x)-\beta T_{j})\\
\geq& \exp(-2|\beta|\|\psi\|)\sum\limits_{j\in \mathbb{N},\
T_{j}-\|\psi\|
>T}\exp ((\beta+\frac{\delta}{2})T_{j}-\beta T_{j})\\=&\infty.
 \end{flalign*}
 Therefore,  for all
$\beta<P_{\psi}(\varphi)-\delta$, \begin{equation}\label{tag-1}
\lim \limits_{\epsilon\rightarrow 0}\limsup \limits_{T \rightarrow
\infty}R_{\psi ,T}(f,\varphi,\{-\beta \|\psi\| m_n\}_{n\in G_T},
\epsilon)=\infty.
\end{equation}
This argument is not only valid for $P_{\psi}(\varphi)\in
\mathbb{R}$, but also for $P_{\psi}(\varphi)=\infty$, in which
case (3.7) holds for every $\beta\in \mathbb{R}$. Then
\begin{equation}\label{tag-1}
P_{\psi}(\varphi)\leq \inf\{\beta \in \mathbb{R}: \lim
\limits_{\epsilon\rightarrow 0}\limsup \limits_{T \rightarrow
\infty}R_{\psi ,T}(f,\varphi,\{-\beta \|\psi\| m_n\}_{n\in G_T},
\epsilon)<\infty\}.
\end{equation}

 Next, we establish the reverse inequality.  We consider the case   $P_{\psi}(\varphi)\in \mathbb{R}$ and show that for any $\delta>0,$
 $$\lim
\limits_{\epsilon\rightarrow 0}\limsup \limits_{T \rightarrow
\infty}R_{\psi ,T}(f,\varphi,\{-(P_{\psi}(\varphi)+\delta)
\|\psi\| m_n\}_{n\in G_T}, \epsilon)<\infty.$$ Again, by the
equivalent definition of $P_{\psi}(\varphi)$,  we have, for any
$\epsilon>0$,
$$\limsup \limits_{T \rightarrow \infty}\frac{1}{T}\log P_{\psi
,T}(f,\varphi, \epsilon)< P_{\psi}(\varphi)+\frac{\delta}{2},$$
and we can find  an $l_{0}\in \mathbb{N}$ such that for all $l \in
\mathbb{N}$ with $l \geq l_{0}$,
$$P_{\psi,lm}(f,\varphi,\epsilon)\leq
\exp(lm(P_{\psi}(\varphi)+\frac{2\delta}{3})).$$
  Note   that for $n\in
S_{lm}, x\in E_{n}$, we have
$$|\|\psi\|
m_n(x)-lm|<2\|\psi\|$$ and
$$-(P_{\psi}(\varphi)+\delta)\|\psi\| m_n(x)\leq -lm(P_{\psi}(\varphi)+\delta)+2|P_{\psi}(\varphi)+\delta||\psi\|.$$
 Moreover, for sufficiently large $T>0, n\in G_{T}, x\in E^{'}_{n} \subset
Y_n$, there exists a unique $l\in \mathbb{N}$ such that
$(l-1)m<S_{n}\psi(x)\leq lm$. Obviously $S_{n+1}\psi(x)>lm$.
Hence, we obtain
\begin{flalign*}
&R_{\psi,T}(f,\varphi,\{-(P_{\psi}(\varphi)+\delta)\|\psi\|
m_n\}_{n\in G_T},\epsilon)\\
\leq& \sum\limits_{l\geq l_{0}}\sup\Bigg\{\sum\limits_{n\in
S_{lm}}\sum \limits_{x\in E_{n}}\exp ((S_{n}\varphi)(x)-
(P_{\psi}(\varphi)+\delta)\|\psi\| m_n(x)):\\
&\quad\quad\quad\quad\quad\quad\quad\quad\quad\quad\quad\quad\quad\quad\quad\quad\quad
E_n \text{ is an } (n,\epsilon)\text{-separated set of } X_n, n\in S_{lm}\Bigg\}\\
\leq& \exp(2\|\psi\||P_{\psi}(\varphi)+\delta|)\sum\limits_{l\geq
l_{0}}\exp(- (P_{\psi}(\varphi)+\delta)lm)P_{\psi,lm}(f,\varphi,\epsilon)\\
\leq& \exp(2\|\psi\||P_{\psi}(\varphi)+\delta|)\sum\limits_{l\geq
l_{0}}\exp(-\frac{\delta}{3}lm)\\<&\exp(2\|\psi\||P_{\psi}(\varphi)+\delta|)\frac{1}{1-\exp
(-\frac{\delta m}{3})}.
\end{flalign*}
 This
implies
$$\lim \limits_{\epsilon\rightarrow 0}\limsup \limits_{T
\rightarrow \infty}R_{\psi
,T}(f,\varphi,\{-(P_{\psi}(\varphi))+\delta) \|\psi\| m_n\}_{n\in
G_T}, \epsilon)<\infty,$$ and hence,
 \begin{equation}\label{tag-1}
 P_{\psi}(\varphi)\geq
\inf\{\beta \in \mathbb{R}: \lim \limits_{\epsilon\rightarrow
0}\limsup \limits_{T \rightarrow \infty}R_{\psi
,T}(f,\varphi,\{-\beta \|\psi\| m_n\}_{n\in G_T},
\epsilon)<\infty\}.
\end{equation} Combining (3.8) and (3.9) we obtain
(3.6).
\begin {cor}
Let $(X,f)$ be a TDS, and $\varphi, \psi \in C(X,\mathbb{R})$ with
$\psi>0$. We have
\begin{equation}\label{tag-1}
P_{\psi}(\varphi)\geq \inf \{\beta \in \mathbb{R}: P(\varphi-\beta
\psi)\leq 0\}.
\end{equation}
\end {cor}
\noindent\textbf{Proof.} Let  $\beta\in \{\beta \in
\mathbb{R}:\lim \limits_{\epsilon\rightarrow 0}\limsup \limits_{T
\rightarrow \infty}R_{\psi ,T}(f,\varphi-\beta \psi,
\epsilon)<\infty\}$ and
$$\lim \limits_{\epsilon\rightarrow 0}\limsup \limits_{T
\rightarrow \infty}R_{\psi ,T}(f,\varphi-\beta \psi,
\epsilon)=a.$$ Then for any   $\epsilon>0$,
$$\limsup \limits_{T \rightarrow \infty}R_{\psi
,T}(f,\varphi-\beta \psi, \epsilon)<a+1 .$$ We  can find a $T_{0}>0$
 such that for all $T>T_{0}$,
$$R_{\psi ,T}(f,\varphi-\beta\psi, \epsilon)<a+2.$$ Now, for
sufficiently large $n\in\mathbb{N}$,
$$S_{n}\psi(x)>T, \ \ \forall x\in X,$$ and hence, for such $n\in G_{T}$, $E_{n}$ is an $(n, \epsilon)$-separated  set of
$X$ and
$$\sum\limits_{x\in E_{n}}\exp (S_{n}(\varphi-\beta \psi))(x)<a+2.$$ It follows
from this that $$P(\varphi-\beta \psi)\leq 0.$$ Since
\begin{flalign*}&\inf \{\beta \in \mathbb{R}:\lim
\limits_{\epsilon\rightarrow 0}\limsup \limits_{T \rightarrow
\infty}R_{\psi ,T}(f,\varphi-\beta \psi, \epsilon)<\infty\} \\
\geq& \inf\{\beta\in\mathbb{R}:P(\varphi-\beta \psi)\leq
0\},\end{flalign*} the inequality (3.10) follows by Theorem 3.1.
\begin {cor}
Let $(X,f)$ be a TDS, and $\varphi, \psi \in C(X,\mathbb{R})$ with
$\psi>0$. We have
$$P_{\psi}(\varphi)=\inf \{\beta \in \mathbb{R}:
P(\varphi-\beta\psi)\leq 0\}=\sup\{\beta \in \mathbb{R}:
P(\varphi-\beta\psi)\geq 0\}.$$
\end {cor}
\noindent\textbf{Proof.} If there exists a $\beta \in \mathbb{R}$
such that $P(\varphi-\beta \psi)=\infty$, then  $P(\varphi-\beta
\psi)=\infty$ for all $\beta \in \mathbb{R}$.  By Corollary 3.1,
we have
$$P_{\psi}(\varphi)=\inf \{\beta \in \mathbb{R}:
P(\varphi-\beta\psi)\leq 0\}=\sup\{\beta \in \mathbb{R}:
P(\varphi-\beta\psi)\geq 0\}.$$ Suppose for  any $\beta \in
\mathbb{R}$, $P(\varphi-\beta \psi)<\infty$.  By (1.1) we have
$$P(\varphi-\beta \psi)=\sup \{h_{\nu}(f)+\int \varphi d\nu- \beta\int\psi d\nu: \nu\in M(X,f)\}.$$ Then for  each
$\beta_{1},\beta_{2} \in \mathbb{R}, \beta_{1}<\beta_{2}$ and
$0<\epsilon<\frac{m(\beta_{2}-\beta_{1})}{2}$, there exists a $\mu
\in M(X,f)$ such that
\begin{flalign*}
&\sup \{h_{\nu}(f)+\int \varphi d\nu- \beta_{2}\int\psi d\nu:
\nu\in M(X,f)\} \\<& h_{\mu}(f)+\int \varphi d\mu-
\beta_{2}\int\psi d \mu +\epsilon\\=&h_{\mu}(f)+\int \varphi d\mu-
\beta_{1}\int\psi d \mu +\epsilon-(\beta_{2}-\beta_{1})\int\psi
d\mu\\<&h_{\mu}(f)+\int \varphi d\mu- \beta_{1}\int\psi d
\mu-(\beta_{2}-\beta_{1})(\int\psi d\mu-\frac{m}{2})\\ \leq &\sup
\{h_{\nu}(f)+\int \varphi d\nu- \beta_{1}\int\psi d\nu: \nu\in
M(X,f)\}-(\beta_{2}-\beta_{1})(\int\psi d\mu-\frac{m}{2}).
\end{flalign*}
Thus, the map $\beta \mapsto P(\varphi-\beta \psi)$ is strictly
decreasing.

Next, we prove that $$P(\varphi-\beta\psi)< 0\Longrightarrow
R_{\psi ,T}(f,\varphi-\beta \psi, \epsilon)<\infty.$$ Let
$P(\varphi-\beta\psi)=2a< 0$. For any $\epsilon>0$, we can find
$N\in \mathbb{N}$  such that for  all
 $n\in \mathbb{N}$ with $n\geq N$,   $$\sup\limits_{E_{n}}\sum
\limits_{x\in E_{n}}\exp (S_{n}(\varphi-\beta\psi))(x)\leq \exp
(na),$$ where the supremum is taken over all
$(n,\epsilon)$-separated sets of $X$. Consequently, for
sufficiently large $T>0$, we have
$$R_{\psi ,T}(f,\varphi-\beta \psi, \epsilon)\leq
\sum\limits_{n\geq N}\sup\limits_{E_{n}}\sum\limits_{x\in
E_{n}}\exp (S_{n}(\varphi-\beta\psi))(x)\leq
\frac{1}{1-\exp(a)}<\infty,$$ and the conclusion holds. \\ Since
$$\inf \{\beta \in \mathbb{R}: P(\varphi-\beta\psi)< 0\}\geq
\inf\{\beta \in \mathbb{R}: \lim \limits_{\epsilon\rightarrow
0}\limsup \limits_{T \rightarrow \infty}R_{\psi
,T}(f,\varphi-\beta \psi, \epsilon)<\infty\},$$ by Theorem 3.1 and
Corollary 3.1, we conclude that
\begin{flalign*}\inf \{\beta \in \mathbb{R}:
P(\varphi-\beta\psi)\leq 0\}&= \inf \{\beta \in \mathbb{R}:
P(\varphi-\beta\psi)< 0\}\\ &=\sup\{\beta \in \mathbb{R}:
P(\varphi-\beta\psi)\geq 0\}.\end{flalign*}

\begin{cor}
Let $(X,f)$ be a TDS and $\varphi, \psi \in C(X,\mathbb{R})$ with
$\psi>0$. Suppose that for each $\beta\in \mathbb{R}$ we have
$P(\varphi-\beta\psi)\in \mathbb{R}$.  Then
$P(\varphi-P_{\psi}(\varphi)\psi)=0$.
\end{cor}
\noindent\textbf{Proof.} By the proof of Corollary 3.2 the map
$\beta \mapsto P(\varphi-\beta \psi)$ is a strictly decreasing,
continuous map on $ \mathbb{R}$. Hence
$P(\varphi-P_{\psi}(\varphi)\psi)=0$.

 We are
now ready to prove Theorem 1.1.\\ \emph{Proof of Theorem 1.1}.
Firstly,  we show
\begin{equation}
P_{\psi}(\varphi)\geq \sup \{\frac{h_{\upsilon}(f)}{\int \psi
d\nu}+\frac{\int \varphi d\nu}{\int \psi d\nu}: \nu \in M(X,f)\}.
\end{equation}
 By Corollary 3.1 we have $0\geq P(\varphi-\beta\psi)$ for $\beta>
 P_{\psi}(\varphi)$.
It follows from (1.1) that
 \begin{flalign*} 0&\geq
P(\varphi-\beta\psi)\\&= \sup \{h_{\nu}(f)+\int \varphi d\nu-
\beta \int\psi d\nu: \nu\in M(X,f)\} \\&=\sup \{\int \psi
d\nu(\frac{h_{\upsilon}(f)}{\int \psi d\nu}+\frac{\int \varphi
d\nu}{\int \psi d\nu}-\beta): \nu \in M(X,f)\},
\end{flalign*}
and hence (3.11) holds.

 Next, we establish the reverse inequality.
 Similarly by Corollary
3.2 we have $ P(\varphi-\beta\psi)\geq 0$ for $
\beta<P_{\psi}(\varphi)$.  Then
\begin{flalign*}
 &P(\varphi-\beta\psi)\\=& \sup \{h_{\nu}(f)+\int \varphi d\nu-
\beta \int\psi d\nu: \nu\in M(X,f)\} \\=&\sup \{\int \psi
d\nu(\frac{h_{\upsilon}(f)}{\int \psi d\nu}+\frac{\int \varphi
d\nu}{\int \psi d\nu}-\beta): \nu \in M(X,f)\}\\
\geq& 0.
\end{flalign*}
It is easy to see that
\begin{equation} P_{\psi}(\varphi)\leq\sup
\{\frac{h_{\upsilon}(f)}{\int \psi d\nu}+\frac{\int \varphi
d\nu}{\int \psi d\nu}: \nu \in M(X,f)\}.
\end{equation} Combining
(3.11) and (3.12), we obtain (1.3).
\section{ A SPECIAL CASE (BS-DIMENSION)}
\quad\quad In this section we will show that the BS dimension
with Carath\'{e}odory structure is a special case of the  induced
pressure. The  BS dimension was first defined by Barreira and
Schmeling [2] as follows.

For $n\geq1, \epsilon>0$, we put
$$\mathcal{W}_{n}(\epsilon)=\{B_{n}(x,\epsilon):x\in X\}.$$
  For
any $B_{n}(x,\epsilon)\in\mathcal{W}_{n}(\epsilon), \psi\in
C(X,\mathbb{R})$ with $\psi>0$, the function $\psi$ can induce a
function by
$$\psi(B)=\sup\limits_{x\in B}(S_{n}\psi)(x).$$
We call $\mathcal{G}\subset \cup_{j\geq
N}\mathcal{\mathcal{W}}_{j}(\epsilon)$  covers $X$, if
$\bigcup\limits_{B\in\mathcal{G}}B=X$.
\begin{defn}
Let $(X,f)$ be a TDS. For any $\alpha>0, N\in\mathbb{N}$ and
$\epsilon>0$, we define $$M(\alpha,\epsilon,N)=\inf
\limits_{\mathcal{G}}\{\sum\limits_{B\in\mathcal{G}}\exp(-\alpha
\psi(B))\},$$ where the infimum is taken over all finite
$\mathcal{G}\subset \cup_{j\geq
N}\mathcal{\mathcal{W}}_{j}(\epsilon)$ that cover $X.$  Obviously
 $M(\alpha,\epsilon,N)$ is a finite outer measure on $X$ and
increases as $N$ increases. Define
$$m(\alpha,\epsilon)=\lim\limits_{N \rightarrow
\infty}M(\alpha,\epsilon,N)$$ and $$\dim_{BS}(X,\epsilon)=\inf
\{\alpha:m(\alpha,\epsilon)=0\}=\sup\{\alpha:m(\alpha,\epsilon)=\infty\}.$$
The BS dimension is $\dim_{BS}X=\lim \limits_{\epsilon\rightarrow
0}\dim_{BS}(X,\epsilon):$ this limit exists because given
$\epsilon_{1}<\epsilon_{2},$ we have $m(\alpha,\epsilon_{1})\geq
m(\alpha,\epsilon_{2})$, so
$\dim_{BS}(X,\epsilon_{1})\geq\dim_{BS}(X,\epsilon_{2})$.
\end {defn}
\begin{prop}
For a TDS, we have $P_{\psi}(0)=\dim_{BS}X$.
\end{prop}
\noindent\textbf{Proof.} By [2, Proposition 6.4], we have
$P(-\psi\dim_{BS}X)=0$. Now it follows from Corollary 3.3 that
$P_{\psi}(0)=\dim_{BS}X$.
\section{EQUILIBRIUM MEASURES AND GIBBS MEASURES  }
\quad\quad In this section we  consider  the problem of the
existence of equilibrium measures for the induced pressure.  We also
study  the relation between Gibbs measures and equilibrium
measures for the induced pressure in the particular case of symbolic
dynamics.
\begin{defn}
Let $(X,f)$ be a TDS and $\varphi,\psi\in C(X,\mathbb{R})$ with
$\psi>0$. A member $\mu$ of $M(X,f)$ is called an equilibrium
measure for $\psi$ and $\varphi$ if
$P_{\psi}(\varphi)=\frac{h_{\mu}(f)+\int \varphi d\mu}{\int \psi
d\mu}.$ We will write $M_{\psi,\varphi}(X,f)$ for the collection
of all equilibrium measures for $\psi$ and $\varphi$.
\end{defn}
\begin{defn}
Let $(X,f)$ be a TDS. Then $f$ is said to be positively expansive
if there exists $\epsilon>0$ such that $x=y$ whenever
$d(f^{n}(x),f^{n}(y))<0$ for every $n\in\mathbb{N}\cup \{0\}$.
\end{defn}
 The entropy map of a TDS is the map $\mu\mapsto h_{\mu}(f)$,
which is defined on $M(X,f)$ and has values in $[0,\infty]$. The
entropy map $\mu\mapsto h_{\mu}(f)$ is called upper
semi-continuous if given a measure $\mu \in M(X,f)$ and
$\delta>0$, we have $h_{\nu}(f)<h_{\mu}(f)+\delta$ for any measure
$\nu\in M(X,f)$ in some open neighborhood of $\mu$. Now we show
that any expansive map has equilibrium measures.
\begin{prop}
Let $(X,f)$ be a TDS and $\varphi,\psi\in C(X,\mathbb{R})$ with
$\psi>0$. Then

{\em ($\romannumeral1$)} If $f$ is a positively expansive map,
then $M_{\psi,\varphi}(X,f)$ is compact and non-empty.

{\em ($\romannumeral2$)} If $\varphi,\phi,\psi\in C(X,\mathbb{R})$
with $\psi>0$ and if there exists a $c\in \mathbb{R}$ such that
$$\varphi-\phi-c\int\psi d\mu\in \overline{\{\tau\circ f
-\tau:\tau\in C(X,\mathbb{R})\}}$$ for each $\mu\in M(X,f)$, then
$M_{\psi,\varphi}(X,f)=M_{\psi,\phi}(X,f)$.
\end{prop}
\noindent\textbf{Proof.} \noindent($\romannumeral1$) For a
positively expansive map $f$,
 it follows from the proof in [1,9] that the
map $\mu \mapsto h_{\mu}(f)$ is upper semi-continuous.  Then
$\mu\mapsto \frac{h_{\mu}(f)}{\int \psi d\mu}$ is upper
semi-continuous. Since the map $$\mu\mapsto \int
\frac{\varphi}{\int\psi d\mu}d\mu$$ is continuous for each
$\varphi\in C(X,f),$ then $$\mu\mapsto
\frac{h_{\mu}(f)+\int\varphi d\mu}{\int\psi d\mu}$$ is upper
semi-continuous. Since an upper semi-continuous map has a maximum
on any compact set, it follows from Theorem 1.1 that
$M_{\psi,\varphi}(X,f)\neq \emptyset$. The upper semi-continuity
also implies $M_{\psi,\varphi}(X,f)$ is compact because if
$\mu_{n}\in M_{\psi,\varphi}(X,f)$ and $\mu_{n}\rightarrow \mu \in
M(X,f)$, then $$\frac{h_{\mu}(f)+\int \varphi d\mu}{\int \psi
d\mu}\geq
\limsup\limits_{n\rightarrow\infty}\frac{h_{\mu_{n}}(f)+\int
\varphi d\mu_{n}}{\int \psi d\mu_{n}}=P_{\psi}(\varphi),$$ so
$\mu\in M_{\psi,\varphi}(X,f)$.
\\
  ($\romannumeral2$) Note
that for each $\mu\in M(X,f)$
$$\frac{h_{\mu}(f)+\int\varphi d\mu}{\int\psi d\mu}=\frac{h_{\mu}(f)+\int\phi d\mu}{\int\psi
d\mu},$$   therefore $P_{\psi}(\varphi)=P_{\psi}(\phi)+c$, hence
$M_{\psi,\varphi}(X,f)=M_{\psi,\phi}(X,f)$.

Next,  we consider symbolic dynamics. Let $(\Sigma_{A} ,\sigma)$
be a one-sided \textit{topological Markov shift}\ (TMS, for short)
over a finite set $S=\{1,2,\ldots, k\}$. This means that there
exists a matrix $A=(t_{ij})_{k\times k}$ of zeros and ones (with
no row or column made entirely of zeros) such that
$$\Sigma_{A} =\{{\omega =(i_{1},i_{2},\ldots)\in
S^{\mathbb{N}}:t_{i_{j}i_{j+1}=1} \text{ for every }
 j\in \mathbb{N}}\}.
$$ The \textit{shift map } $\sigma :\Sigma_{A} \rightarrow\Sigma_{A}$ is
defined by $(i_{1},i_{2},i_{3}\ldots)\mapsto
(i_{2},i_{3},\ldots)$. We call $C_{i_{1}\ldots
i_{n}}=\{(j_{1}j_{2}\ldots) \in \Sigma_{A} :j_{l}=i_{l} \text{ for
}l=1,\ldots,n\}$ the \textit{cylindrical set } of $\omega$. We
equip $\Sigma_{A}$ with the topology generated by the cylindrical
sets. The topology of a TMS is metrizable and may be given by the
metric $d_{\alpha}(\omega
,\omega')=e^{-\alpha|\omega\wedge\omega'|}, \alpha>0$, where
$\omega\wedge\omega'$ denotes the longest common initial block of
$\omega ,\omega'\in \Sigma_{A}$. The shift map $\sigma$ is
continuous with respect to this metric.  A TMS
$(\Sigma_{A},\sigma)$ is called a topologically mixing TMS if for
every $a, b \in S$, there exists an $N_{ab}\in\mathbb{N}$ such
that for every $n>N_{ab}$, we have $C_{a}\cap \sigma^{-n}C_{b}$.
\begin{defn}
Let $(\Sigma_{A} ,\sigma)$ be a TMS and $\varphi,\psi\in
C(\Sigma_{A},\mathbb{R})$ with $\psi>0$.  We say that a probability
measure $\mu$ in $\Sigma_{A}$ is a Gibbs measure for $\psi$ and
$\varphi$
 if there exists a $K>1$ such that $$K^{-1}\leq
\frac{\mu(C_{i_{1}\ldots i_{n}})}{\exp [-(S_{n}\psi)(\omega)
P_{\psi}(\varphi)+(S_{n}\varphi)(\omega)]}\leq K$$ for each
$(i_{1},i_{2},\ldots)\in \Sigma_{A}, n\in \mathbb{N}$ and
$\omega\in C_{i_{1}\ldots i_{n}}$.
\end{defn}
 We show that $\sigma$-invariant Gibbs measures are equilibrium
measures. Making a similar proof as in [1, Theorem 3.4.2], we can
obtain the following statement:
\begin{prop}
If a probability measure $\mu$ in $(\Sigma_{A},\sigma)$ is a
$\sigma$-invariant Gibbs measure for $\varphi$ and $\psi$, then it
is also an equilibrium measure for  $\psi$ and $\varphi$.
\end{prop}  Now we establish the existence of Gibbs measures.
\begin{prop}
Let $(\Sigma_{A},\sigma)$ be a topologically mixing TMS. Suppose
that $\varphi$ and $\psi$ are H$\ddot{\text{o}}$lder continuous
functions and $\psi>0$. Then there exists at least one
$\sigma$-invariant Gibbs measure for  $\psi$ and $\varphi$.
\end{prop}
\noindent\textbf{Proof.} By Corollary 3.3 we have
$$P(\varphi-P_{\psi}(\varphi)\psi)=0.$$ As
$\varphi-P_{\psi}(\varphi)\psi$ is H$\ddot{\text{o}}$lder
continuous, it follows from [1, Theorem 3.4.4] that there
exists a $K>1$ such that $$K^{-1}\leq \frac{\mu(C_{i_{1}\ldots
i_{n}})}{\exp [-n
P(\varphi-P_{\psi}(\varphi)\psi)-(S_{n}\psi)(\omega)
P_{\psi}(\varphi)+(S_{n}\varphi)(\omega)]}\leq K.$$

\noindent {\bf ACKNOWLEDGEMENTS.}  This  research was supported by the
National Natural Science Foundation of China (Grant No.
11271191) and the National Basic Research Program of China (Grant No.
2013CB834100).  We would like to thank the referee for very useful
comments and helpful suggestions.  The first author would like to
thank Dr. Zheng Yin for useful discussions.

\end{document}